\documentclass[11pt]{article}

\usepackage{amssymb,amscd,array}

\let\ssection=\section
\renewcommand{\section}{\setcounter{equation}{0}\ssection}

\newcommand{\bbR}{\mathbb{R}}
\newcommand{\bbRP}{\mathbb{RP}}
\newcommand{\bbCP}{\mathbb{CP}}

\newcommand{\bbC}{\mathbb{C}}
\newcommand{\bbN}{\mathbb{N}}

\newcommand{\cD}{{\mathcal{D}}}

\newcommand{\cE}{{\mathcal{E}}}
\newcommand{\Diff}{\mathrm{Diff}}
\newcommand{\Div}{\mathrm{Div}}
\newcommand{\rE}{\mathrm{E}}

\newcommand{\cF}{{\mathcal{F}}}
\newcommand{\gl}{{\mathrm{gl}}}

\newcommand{\Hom}{\mathrm{Hom}}

\newcommand{\cL}{{\mathcal{L}}}

\newcommand{\Pol}{\mathrm{Pol}}

\newcommand{\cS}{{\mathcal{S}}}

\newcommand{\SL}{\mathrm{SL}}
\newcommand{\Sl}{\mathrm{sl}}

\newcommand{\so}{\mathrm{o}}

\newcommand{\Vect}{\mathrm{Vect}}

\newcommand{\cqfd}{\hspace*{\fill}\rule{3mm}{3mm}}

\begin{document}

%\baselineskip=18pt

%\textwidth=16truecm
%\textheight=24truecm
%\hoffset=-1.5truecm
%\voffset=-2.5truecm

\def\d{\delta}
\def\g{\gamma}
\def\om{\omega}
\def\r{\rho}
\def\a{\alpha}
\def\s{\sigma}
\def\vfi{\varphi}
\def\l{\lambda}
\def\m{\mu}
\def\implies{\Rightarrow}

\oddsidemargin .1truein
\newtheorem{thm}{Theorem}[section]
\newtheorem{lem}[thm]{Lemma}%[section]
\newtheorem{cor}[thm]{Corollary}%[section]
\newtheorem{pro}[thm]{Proposition}%[section]
\newtheorem{ex}[thm]{Example}%[section]
\newtheorem{rmk}[thm]{Remark}%[section]
\newtheorem{defi}[thm]{Definition}%[section]
%\newremark{ex}[thm]{Example}

%\newtheorem{thm}{Theorem}
%\newtheorem{lem}{Lemma}
%\newtheorem{cor}{Corollary}
%\newtheorem{prop}[thm]{Proposition}
%\newtheorem{definition}{Definition}

\title{Projectively equivariant symbol calculus}

\author{P.B.A. Lecomte
\thanks{Institute de Math\'ematiques,
 Universit\'e de Li\`ege, Sart Tilman,
Grande Traverse, 12 (B 37),
B-4000 Li\`ege,
BELGIUM, mailto:plecomte@ulg.ac.be
}
\and
 V.Yu. Ovsienko
\thanks{C.N.R.S., Centre de Physique Th\'eorique,
 Luminy -- Case 907,
F--13288 Marseille, Cedex 9,
FRANCE, mailto:ovsienko@cpt.univ-mrs.fr
}
}

\date{}

\maketitle

\thispagestyle{empty}

\begin{abstract}
The spaces of linear differential operators $\cD_\l(\bbR^n)$ acting on
$\l$-densities on $\bbR^n$ and the space 
$\Pol(T^*\bbR^n)$ of functions on $T^*\bbR^n$ 
which are polynomial on the fibers are not isomorphic 
as modules over the Lie algebra $\Vect(\bbR^n)$ 
of vector fields of $\bbR^n$. However, these modules 
are isomorphic as $\Sl(n+1,\bbR)$-modules where
$\Sl(n+1,\bbR)\subset \Vect(\bbR^n)$
is the Lie algebra of infinitesimal projective transformations. In
addition, such an $\Sl_{n+1}$-equivariant bijection 
is unique (up to normalization). This leads to a notion
of projectively equivariant quantization and symbol calculus
for a manifold endowed with a (flat) projective structure.

We apply the $\Sl_{n+1}$-equivariant symbol map to study
the $\Vect(M)$-modules $\cD_\l^k(M)$ 
of $k$-th-order linear differential operators acting on
$\l$-densities, for an arbitrary manifold $M$
and classify the quotient-modules 
$\cD_\l^k(M)/\cD_\l^\ell(M)$.
\end{abstract}

%\newpage

%%%%%%%%%%%%%%%%%%%%%%%%%%%%%%%%%%%%%%%%%%%%%%%%%%%%%%%%%%%%%%%%%%%%%%%%%%%%%%%%%%%%
%%%%%%%%%%%%%%%%%%%%%%%%%%%%%%%%%%%%%%%%%%%%%%%%%%%%%%%%%%%%%%%%%%%%%%%%%%%%%%%%%%%%
\section{Introduction}
%%%%%%%%%%%%%%%%%%%%%%%%%%%%%%%%%%%%%%%%%%%%%%%%%%%%%%%%%%%%%%%%%%%%%%%%%%%%%%%%%%%%
%%%%%%%%%%%%%%%%%%%%%%%%%%%%%%%%%%%%%%%%%%%%%%%%%%%%%%%%%%%%%%%%%%%%%%%%%%%%%%%%%%%%

Roughly speaking, a quantization procedure associates linear
operators on a Hilbert space to functions on a symplectic manifold. In
particular, if this manifold is the cotangent bundle $T^*M$ of a smooth
manifold $M$, a now standard quantization procedure leads to a linear bijection
from the space $\Pol(T^*M)$ of functions on $T^*M$ which are polynomial on the
fibers (or, equivalently, symmetric contravariant tensor fields on $M$),
into the space $\cD(M)$ of linear differential operators on $M$~:
$$
Q :\Pol(T^*M)\to\cD(M)\,,
\qquad
\hbox{quantization map}
$$
The inverse
$$
\sigma = Q^{-1}:\cD(M)\to\Pol(T^*M)\,,
\qquad
\hbox{symbol map}
$$
associates to each operator $A\in\cD(M)$ a sort of {\it total} symbol.

\medskip
\noindent
{\bf 1.1}
One of the main questions usually arising in this context is to find the
group of symmetries, that is, the Lie group acting on $M$ such that
the quantization procedure is {\it equivariant} with respect to this
action. It is natural then to consider the space of differential
operators $\cD(M)$ as a module over the group $\Diff(M)$ of diffeomorphisms
of $M$.

There exists a natural family of $\Diff(M)$-
and $\Vect(M)$-module structures on $\cD(M)$.
To define it, one considers differential operators as acting on
tensor-densities of arbitrary degree $\l$, this leads to a $\Diff(M)$-module
 of differential operators $\cD_\l(M)$.
The module $\cD_{1/2}(M)$ on half-densities
plays a special role (see \cite{bla,kos}).

It worth noticing that modules of differential operators on tensor densities
have been studied in recent papers \cite{duv,lec,gar,cmz}.

\medskip
\noindent
{\bf 1.2}
One of the difficulties of quantization is that there is 
{\it no natural quantization map.} In other words, $\Pol(T^*M)$ and
$\cD_\l(M)$ {\it are not isomorphic} as modules over the group $\Diff(M)$ of
diffeomorphisms of $M$ nor as modules over its
Lie algebra $\Vect(M)$ of vector fields.

The main idea of this paper is to fix a {\it maximal group of
symmetries} for which such quantization is still possible.
We consider a manifold $M$ (of dimension~$n$)
endowed with a (flat) projective structures (i.e.
we locally identify $M$ with the projective space $\bbRP^n$).
We show that there exist {\it unique} (up to a natural normalization)
symbol and quantization maps equivariant with respect to the group
$\SL(n+1,\bbR)$ of projective symmetries.

In terms of modules of differential operators,
our main result can be formulated as follows.
We consider a natural embedding $\Sl(n+1,\bbR)\subset\Vect(\bbR^n)$ 
($\Sl(n+1,\bbR)$ acts on $\bbR^n$ by infinitesimal
projective transformations). It turns out that 
$\cD_\l(\bbR^n)$ {\it and}
$\Pol(T^*\bbR^n)$ {\it are equivalent $\Sl(n+1,\bbR)$-modules}.
In particular, the $\Sl(n+1,\bbR)$-modules $\cD_\l(\bbR^n)$
with different values of $\l$ are isomorphic to each other.

\medskip
\noindent
{\bf 1.3}
We apply our $\Sl_{n+1}$-equivariant symbol map to the problem of
classification of modules of differential operators on an arbitrary
smooth manifold $M$.
The classification of $\Vect(M)$-modules $\cD^k_\l(M)$ 
has been performed in a series of recent works
\cite{duv,lec} and \cite{gar}. In this paper we classify
the quotient-modules $\cD_\l^k(M)/\cD_\l^\ell(M)$.
We prove that for $k-\ell\geq2$ every such module 
(except the special case $\l=1/2$ with \allowbreak $k-\ell=2$) is
a \textit{nontrivial deformation} of the
$\Vect(M)$-module of symbols and compute the 
corresponding cohomology classes of $\Vect(M)$.

\medskip
\noindent
{\bf 1.4}
An important aspect discussed in this paper is the property of
\textit{locality}. We show that any linear map on $\Pol(T^*\bbR^n)$
is necessarily local if it is
equivariant with respect to the Lie group generated by
translations and homotheties of $\bbR^n$
(i.e., with respect to $\bbR^*\ltimes\bbR^n$-action).
This proves, in particular, that an $\Sl_{n+1}$-equivariant symbol
is given by a differential map.

\medskip
\noindent
{\bf 1.5 Remarks}. 

(a) The relationship between  differential
operators and projective geometry has already been studied in the
fundamental book \cite{wil}. The best known example is the
Sturm-Liouville operator $d^2/dx^2+u(x)$ describing a projective structure
on $\bf R$ (or on $S^1$ if $u(x)$ is periodic).

(b) In the one-dimensional case ($n=1$),
the $\Sl_2$-equivariant symbol map and quantization map
were obtained (in a more general situation of pseudodifferential operators)
in recent work by P.B. Cohen, Yu. I. Manin and D. Zagier \cite{cmz}.
If $n=1$, our isomorphisms coincide with those of \cite{cmz}.
This article is a revised version of the electronic preprint
\cite{LO}; we were not aware of the article \cite{cmz} while the computation
of the $\Sl_{n+1}$-equivariant symbol was carried out.

(c) In the (algebraic) case of global differential operators
on $\bbCP^n$, existence and uniqueness of the
$\Sl_{n+1}$-equivariant symbol is a corollary of
Borho-Brylinski's results \cite{bor}:
${\cal D}({\bf CP}^n)$, as a module over $\Sl(n+1,\bbC)$,
has a decomposition as a sum of irreducible submodules of
{\it multiplicity one}. This implies the uniqueness result.
Our explicit formul{\ae} are valid in the holomorphic case and
define an isomorphism between $\cD(\bbCP^n)$
and the space of functions on $T^*\bbCP^n$ which are polynomial on fibers.

\medskip

We believe that the appearance of projective symmetries
in the context of quantization is natural.
We do not know any work on this subject except the special
one-dimensional case.
A particular role of the Lie algebra of projective transformations
is related to the fact that $\Sl(n+1,\bbR)\subset\Vect(\bbR^n)$ is
a \textit{maximal} Lie subalgebra. Any bigger Lie subalgebra of
$\Vect(\bbR^n)$ is infinite-dimensional.

%%%%%%%%%%%%%%%%%%%%%%%%%%%%%%%%%%%%%%%%%%%%%%%%%%%%%%%%%%%%%%%%%%%%%%%%%%%%%%%%%%%%
%%%%%%%%%%%%%%%%%%%%%%%%%%%%%%%%%%%%%%%%%%%%%%%%%%%%%%%%%%%%%%%%%%%%%%%%%%%%%%%%%%%%
\section{Modules of differential operators on $\bbR^n$}
%%%%%%%%%%%%%%%%%%%%%%%%%%%%%%%%%%%%%%%%%%%%%%%%%%%%%%%%%%%%%%%%%%%%%%%%%%%%%%%%%%%%
%%%%%%%%%%%%%%%%%%%%%%%%%%%%%%%%%%%%%%%%%%%%%%%%%%%%%%%%%%%%%%%%%%%%%%%%%%%%%%%%%%%%

%%%%%%%%%%%%%%%%%%%%%%%%%%%%%%%%%%%%%%%%%%%%%%%%%%%%%%%%%%%%%%%%%%%%%%%%%%%%%%%%%%%%
\subsection{Definition of the $\Vect(\bbR^n)$-module structures}\label{module}
%%%%%%%%%%%%%%%%%%%%%%%%%%%%%%%%%%%%%%%%%%%%%%%%%%%%%%%%%%%%%%%%%%%%%%%%%%%%%%%%%%%%

Let us recall the definition of the natural 1-parameter family of
$\Vect(\bbR^n)$-module structures on the
space of differential operators (see \cite{duv,lec,gar}).

\medskip
\noindent
{\bf Definition.}
For each $\l\in\bbR$ (or $\bbC$), consider the space $\cF_\l$ of 
{\it tensor densities}  of degree  $\l$ on $\bbR^n$, that is,
of sections of the line bundle $|\Lambda^n{}T^*{}\bbR^n|^{\otimes\l}$.
The space $\cF_\l$ has a natural structure of $\Vect(\bbR^n)$-module,
defined by the Lie derivative.

\medskip

In coordinates~:
$$
\phi
=
\phi(x^1,\ldots,x^n)\left|dx^1\wedge\ldots\wedge dx^n\right|^\l.
$$
The action of 
$X\in\Vect(\bbR^n )$ on $\phi\in C^\infty(\bbR^n )$ is given by
\begin{equation}
L_X^\l(\phi) = 
X^i \partial_i\phi +
\l\partial_iX^i\phi\,,
\label{LieDer}
\end{equation}
where $\partial_i=\partial/\partial x^i$.
Note, that the formula (\ref{LieDer})
does not depend on the choice of local coordinates.

\medskip
\noindent
{\bf Remark.}
Modules $\cF_\l$ are not isomorphic to each other for different
values of $\l$ (cf. \cite{fuc}).
The simplest examples of modules of tensor densities are 
$\cF_0=C^{\infty}(\bbR^n)$ and $\cF_1=\Omega^n(\bbR^n)$, the module $\cF_{1/2}$
is particularly important for geometric quantization 
(see \cite{bla,kos}).

\medskip
\noindent
{\bf Definition.}
Consider the space $\cD_\l(\bbR^n)$ (or $\cD_\l$ for short)
of differential operators on tensor densities,
$A:\cF_\l\to\cF_\l$.
The natural $\Vect(\bbR^n)$-action on $\cD_\l$ is given by
\begin{equation} 
\label{E6}
\cL_X^\l(A) = 
L_X^\l\circ A - A \circ L_X^\l.
\end{equation}
Denote $\cD^k_\l\subset\cD_\l$ the
$\Vect(\bbR^n)$-module of $k$-th-order differential operators.

\medskip

In local coordinates any linear differential
operator order $k$ is of the form~:
\begin{equation}
\label{E4}
A=a_k^{{i_1}\ldots{i_k}}
\partial_{i_1}\cdots
\partial_{i_k}
+\cdots+
a_1^i\partial_i
+a_0
\end{equation}
with coefficients
$a_j^{{i_1}\ldots{i_j}}=a_j^{{i_1}\ldots{i_j}} (x^1,\dots,x^n)\in
C^{\infty}(\bbR^n)$ (sum over repeated indices is understood).

%%%%%%%%%%%%%%%%%%%%%%%%%%%%%%%%%%%%%%%%%%%%%%%%%%%%%%%%%%%%%%%%%%%%%%%%%%%%%%%%%%%%
\subsection{$\Vect(\bbR^n)$-modules of symmetric tensor fields on $\bbR^n$}
%%%%%%%%%%%%%%%%%%%%%%%%%%%%%%%%%%%%%%%%%%%%%%%%%%%%%%%%%%%%%%%%%%%%%%%%%%%%%%%%%%%%

Consider the space $\Pol(T^*\bbR^n)$ of functions on ${T^*}\bbR^n\cong\bbR^{2n}$
polynomial on the fibers:
$
P=\sum_{j=0}^{k}a_j^{{i_1}\ldots{i_j}}\xi_{i_1}\cdots\xi_{i_j}\,.
$
This space has a natural $\Vect(\bbR^n)$-module structure (defined by the
lift of a vector field to the cotangent space)~:
\begin{equation}
\label{StandAct}
L_X=
X^i\partial_i
-
\xi_j\partial_iX^j\partial_{\xi_i}\,,
\end{equation}
here and below we denote $\partial_i=\partial/\partial x^i$ and
$\partial_{\xi_i}=\partial/\partial \xi_i$.

As a $\Vect(\bbR^n )$-module, $\Pol(T^*{\bbR}^n)$ is
isomorphic to the direct sum of the modules of symmetric contravariant tensor
fields on $\bbR$, i.e.
$
\Pol^k (T^*{\bbR}^n) \cong
{\cal S}^0 \oplus \cdots \oplus {\cal S}^k,
$
where ${\cal S}^\ell = \Gamma(S^\ell(T{\bbR}^n ))$.

%%%%%%%%%%%%%%%%%%%%%%%%%%%%%%%%%%%%%%%%%%%%%%%%%%%%%%%%%%%%%%%%%%%%%%%%%%%%%%%%%%%%
\subsection{Identification of the vector spaces
$\cD_\l$ and $\Pol({T^*}\bbR^n)$}
%%%%%%%%%%%%%%%%%%%%%%%%%%%%%%%%%%%%%%%%%%%%%%%%%%%%%%%%%%%%%%%%%%%%%%%%%%%%%%%%%%%%

\medskip
\noindent
{\bf Definition.}
Substituting the monomial $\xi_{i_1}\cdots\xi_{i_j}$ (where
$\xi=(\xi_i)\in{\bbR^n}^*$) to the partial derivation
$\frac{\partial}{\partial x^{i_1}}\dots\frac{\partial}{\partial x^{i_j}}$,
allows us to identify $A$ with some element of $\Pol(T^*\bbR^n)$~:
\begin{eqnarray}
\label{Identific}
a_k^{{i_1}\ldots{i_k}}
\partial_{i_1}\cdots
\partial_{i_k}
\mapsto
a_k^{{i_1}\ldots{i_k}}
\xi_{i_1}\cdots
\xi_{i_k}
\end{eqnarray}
We get in this way an {\it isomorphism of vector spaces}
$\cD\cong\Pol(T^*\bbR^n)$. 

\medskip

The $\Vect(\bbR^n)$-action (\ref{E6}) is, of course, different from
the standard $\Vect(\bbR^n)$-action (\ref{StandAct}) on polynomials.
We will, therefore, distinguish two $\Vect(\bbR^n)$-modules~:
\begin{eqnarray}
\label{DvaMod}
\cD_\l
&\equiv&
(\Pol(T^*\bbR^n),\cL^\l),\\
\cS
&\equiv&
(\Pol(T^*\bbR^n),L).
\end{eqnarray}

\medskip

In particular, a vector field $X$ corresponds to a first-order polynomial:
$X=X^i\xi_i$. The operator of Lie derivative is then given by a Hamiltonian
vector field
\begin{equation}
L_X=
\partial_{\xi_i}X\partial_i
-
\partial_iX\partial_{\xi_i}\,.
\label{LieDer1}
\end{equation}

\medskip
\noindent
{\bf Remark.}
The identification (\ref{Identific}) is often called in mathematical physics
the {\it normal ordering}. Another frequently used way to identify the spaces
of differential operators on $\bbR^n$ with the space $\Pol(T^*\bbR^n)$ is
provided by the Weyl symbol calculus.

%%%%%%%%%%%%%%%%%%%%%%%%%%%%%%%%%%%%%%%%%%%%%%%%%%%%%%%%%%%%%%%%%%%%%%%%%%%%%%%%%%%%
\subsection{Comparison of the $\Vect(\bbR^n)$-action
on $\cD_\l$ and $\cS$}
%%%%%%%%%%%%%%%%%%%%%%%%%%%%%%%%%%%%%%%%%%%%%%%%%%%%%%%%%%%%%%%%%%%%%%%%%%%%%%%%%%%%

Let us compare the
$\Vect(\bbR^n)$-action on $\cD_\l$ with the standard
$\Vect(\bbR^n)$-action (\ref{LieDer1})
on $\Pol({T^*}\bbR^n)$. We will use the identification (\ref{Identific}) and
write the $\Vect(\bbR^n)$-action (\ref{E6}) in terms of polynomials.

\begin{lem}
The $\Vect(\bbR^n)$-action on $\cD_\l$ has the following form~:
\begin{equation}
\begin{array}{ccl}
\cL_X^\l&=&
L_X-\frac{1}{2}\,\partial_{ij}X
\partial_{\xi_i\xi_j}-
\l\,(\partial_i\circ\Div)X
\partial_{\xi_i}\\[8pt]
&&+(\hbox{higher order derivatives
$\partial_{i_1}\cdots\partial_{i_l}X$})\,,
\end{array}
\label{compar}
\end{equation}
where $\Div\, X=\partial_i\partial_{\xi_i}X=\partial_iX^i$.
\label{PolLem}
\end{lem}
\noindent
{\it Proof}. Straightforward computation.
\cqfd

%%%%%%%%%%%%%%%%%%%%%%%%%%%%%%%%%%%%%%%%%%%%%%%%%%%%%%%%%%%%%%%%%%%%%%%%%%%%%%%%%%%%
%%%%%%%%%%%%%%%%%%%%%%%%%%%%%%%%%%%%%%%%%%%%%%%%%%%%%%%%%%%%%%%%%%%%%%%%%%%%%%%%%%%%
\section{Projective Lie algebra $\Sl(n+1,\bbR)$}
%%%%%%%%%%%%%%%%%%%%%%%%%%%%%%%%%%%%%%%%%%%%%%%%%%%%%%%%%%%%%%%%%%%%%%%%%%%%%%%%%%%%
%%%%%%%%%%%%%%%%%%%%%%%%%%%%%%%%%%%%%%%%%%%%%%%%%%%%%%%%%%%%%%%%%%%%%%%%%%%%%%%%%%%%

The standard action of the Lie algebra $\Sl(n+1,\bbR)$ on $\bbR^n$
is generated by the vector fields~:
\begin{equation} 
\label{E1}
\partial_i\,, 
\quad 
x^i\partial_j\,, 
\quad
x^i{\cal E}\,,
\end{equation}
where
\begin{equation} 
{\cal E} = x^j \partial_j\,.
\label{EuField}
\end{equation}
It will be called the
{\it projective Lie algebra}.
It contains the {\it affine Lie
algebra}, the Lie subalgebra generated by the constant and the linear vector
fields.

\medskip
\noindent
{\bf Remark}.
The group $\SL(n+1,\bbR)$ acts on ${\bbRP}^n$ by homographies
(linear-fractional tranformations). This action is locally defined on
${\bbR}^n \subset {\bbRP}^n$. The Lie subalgebra of $\Vect({\bbR}^n )$
tangent to this action coincides with (\ref{E1}).

%%%%%%%%%%%%%%%%%%%%%%%%%%%%%%%%%%%%%%%%%%%%%%%%%%%%%%%%%%%%%%%%%%%%%%%%%%%%%%%%%%%%
\subsection{Maximality of $\Sl(n+1,\bbR)$ inside $\Vect(\bbR^n)$}
%%%%%%%%%%%%%%%%%%%%%%%%%%%%%%%%%%%%%%%%%%%%%%%%%%%%%%%%%%%%%%%%%%%%%%%%%%%%%%%%%%%%

It seems to be quite a known fact
that the projective Lie algebra $\Sl(n+1,\bbR)$ is 
a maximal subalgebra of the  Lie algebra of 
polynomial vector fields on $\bbR^n$ (cf. \cite{GKO} for $n=1,2$.)
For the sake of completeness, we will prove here the
following version of this result, supposing for simplicity
that $n\geq3$.

\begin{pro}
Given an arbitrary polynomial vector field
$X\not\in\Sl(n+1,\bbR)$, the  Lie algebra generated by
$\Sl(n+1,\bbR)$ and $X$ is
the  Lie algebra of all polynomial vector fields on $\bbR^n$.
\label{maximal}
\end{pro}

\noindent
\textit{Proof.} We will need the following lemmas.
\begin{lem}
Given a polynomial vector field $X$ such that for every $i=1,\ldots,n$
one has: $[\partial_i,X]=\partial_iX\in\Sl(n+1,\bbR)$,
each component of the vector field $X$ is a polynomial in $x^1,\ldots,x^n$
of degree at most 2.
\label{maxlem1}
\end{lem}
\noindent
\textit{Proof.}
Suppose that $\deg X\geq3$ and that $\partial_iX=\a_i\cE$, where $\a_i$ are
some linear functions and $\cE$ is the Euler field~(\ref{EuField}).
From $\partial_i\partial_j=\partial_j\partial_i$, one has:
$(\partial_i\a_j-\partial_i\a_j)\cE+\a_j\partial_i-\a_i\partial_j=0$.
If $n\geq3$, this system has a unique solution $\a_1=\cdots=\a_n\equiv0$.
\cqfd

\medskip

We will also need the following classical fact~:

\begin{lem}
The space $\cS^1_2$ of symmetric 1-contravariant, 2-covariant tensor
field on $\bbR^n$
is split as an $\gl(n,\bbR)$-module into two irreducible components:
$$
\cS^1_2={\cal A}\oplus{\cal B}\,,
\qquad
\hbox{where}
\quad 
{\cal A}=\bbR^n{}^*\cdot\cE\subset\Sl(n+1,\bbR)
\quad
\hbox{and}
\quad
{\cal B}={\rm Ker}(\Div)\,,
$$
and where $\Div=\partial_i\partial_{\xi_i}$.
\label{maxlem2}
\end{lem}
\noindent

\medskip

Let us now prove the proposition.
Considering the commutators $[\partial_i,X]$ and using Lemma \ref{maxlem1},
we can suppose that the vector field $X$ is a homogeneous polynomial
of degree 2. Then, $X=X_A+X_B$ where $X_A\in{\cal A}$ and
$X_B\in{\cal B}$. By Lemma \ref{maxlem2}, the Lie algebra
$\Sl(n+1,\bbR)$ and $X_B$ generate the space of all second order
vector fields on $\bbR^n$.
Proposition \ref{maximal} follows then from the following lemma
whose proof is straightforward.

\begin{lem}
If $n\geq2$, then the space of all vector fields on $\bbR^n$ polynomial
of degree 2 generates the Lie algebra of all polynomial vector fields
on $\bbR^n$.
\label{maxlem3}
\end{lem}

Proposition \ref{maximal} is proven.
\cqfd

\begin{cor}
\hfill\break
(i)
The projective Lie algebra $\Sl(n+1,\bbR)$ is 
a maximal subalgebra of the Lie algebra of polynomial vector fields.
\hfill\break
(ii)
The projective Lie algebra $\Sl(n+1,\bbR)$ is 
a maximal subalgebra of $\Vect(\bbR^n)$ in the class of
finite-dimensional Lie subalgebras.
\label{maximal}
\end{cor}

\medskip
\noindent
\textbf{Remark}.
It worth noticing that
the problem of classification of finite-dimensional Lie subalgebras of
$\Vect(\bbR^n)$ goes back to S.~Lie and remains open (cf. \cite{GKO}).
The only known examples of maximal semisimple Lie subalgebras of
$\Vect(\bbR^n)$ are the projective Lie algebra $\Sl(n+1,\bbR)$ and
the conformal Lie algebra $\so(p+1,q+1)$ with $p+q=n$.
The second case was studied in \cite{DO} as a continuation of \cite{LO}.

%%%%%%%%%%%%%%%%%%%%%%%%%%%%%%%%%%%%%%%%%%%%%%%%%%%%%%%%%%%%%%%%%%%%%%%%%%%%%%%%%%%%
\subsection{Action of $\Sl(n+1,\bbR)$ on $\cD_\l$}
%%%%%%%%%%%%%%%%%%%%%%%%%%%%%%%%%%%%%%%%%%%%%%%%%%%%%%%%%%%%%%%%%%%%%%%%%%%%%%%%%%%%

The action (\ref{E6}) restricted to $\Sl(n+1,\bbR)$ 
takes a particularly nice form.

\begin{pro}
The action of $\Sl(n+1,\bbR)$ on $\cD_\l$
reads~:
\begin{eqnarray}
\cL_X^\l&=&L_X,\label{first}\\
\cL_{X_s}^\l&=&
L_{X_s}-\Big(\rE+\l(n+1)\Big)\partial_{\xi_s}
\label{SLAction}
\end{eqnarray}
for $X$ from the affine subalgebra of $\Sl(n+1,\bbR)$ and
$X_s$ a quadratic vector field given by (\ref{E1})
and where
\begin{equation}
\label{rE}
\displaystyle
\rE=
\xi_i\partial_{\xi_i}\,.
\end{equation}
\label{SLPro}
\end{pro}

\noindent
\textit{Proof.} 
Straightforward in view of Lemma \ref{PolLem}.
\cqfd

\medskip

The formula (\ref{first}) implies that the identification (\ref{Identific})
is an isomorphism between $\Pol({T^*}\bbR^n)$ and $\cD_\l$ as modules
over the affine Lie algebra.

%%%%%%%%%%%%%%%%%%%%%%%%%%%%%%%%%%%%%%%%%%%%%%%%%%%%%%%%%%%%%%%%%%%%%%%%%%%%%%%%%%%%
%%%%%%%%%%%%%%%%%%%%%%%%%%%%%%%%%%%%%%%%%%%%%%%%%%%%%%%%%%%%%%%%%%%%%%%%%%%%%%%%%%%%
\section{Projectively equivariant symbol map}
%%%%%%%%%%%%%%%%%%%%%%%%%%%%%%%%%%%%%%%%%%%%%%%%%%%%%%%%%%%%%%%%%%%%%%%%%%%%%%%%%%%%
%%%%%%%%%%%%%%%%%%%%%%%%%%%%%%%%%%%%%%%%%%%%%%%%%%%%%%%%%%%%%%%%%%%%%%%%%%%%%%%%%%%%

\goodbreak

\medskip
\noindent
{\bf Definition}.

(a)
A linear bijection
$
\s:\cD_\l\to\Pol({T^*\bbR}^n)
$
is called a {\it symbol map} if for every $A\in\cD_\l$,
the highest order term of $\s(A)$ coincides with 
the \textit{principal symbol} 
\begin{equation}
\label{Princip}
\s_A
=
a_k^{i_1 \ldots i_k}\xi_{i_1}\cdots\xi_{i_k}.
\end{equation}

(b)
We will say that $\s$ is
\textit{differential} if in addition its
restriction to each homogeneous component of $\cD_\l$
is a differential operator.

\medskip
\noindent
{\bf Remark}.
The coefficients $a_j^{i_1 \ldots i_j}$ of the differential operator
(\ref{E4}) have no geometric meaning, except for $j=k$.
In other words, there is no natural symbol map, i.e.
there is no linear bijection from
$\cD_\l$ into the space of symmetric contravariant tensor fields on $\bbR^n$
equivariant with respect to the $\Vect(\bbR^n)$-action.

We will see below that for $\l=1/2$, there is also a well-defined
symbol of order $k-1$ on $\cD_{1/2}^k$.

%%%%%%%%%%%%%%%%%%%%%%%%%%%%%%%%%%%%%%%%%%%%%%%%%%%%%%%%%%%%%%%%%%%%%%%%%%%%%%%%%%%%
\subsection{The divergence operator as an affine invariant}
%%%%%%%%%%%%%%%%%%%%%%%%%%%%%%%%%%%%%%%%%%%%%%%%%%%%%%%%%%%%%%%%%%%%%%%%%%%%%%%%%%%%

Let us introduce a differential operator
on the space of polynomials $\cS\cong\Pol({T^*}\bbR^n)$
(and on $\cD_\l$ using the identification (\ref{Identific}))~:
\begin{equation}
\Div=
\partial_i\partial_{\xi_i}\,.
\label{diver}
\end{equation}

On a homogeneous component of order $k$, $\cS^k\subset\cS$ one has~:
%$\Div=k\,{\mathrm{div}}$, where
\begin{equation}
\Div\Big|_{\cS^k}
=
k\,{\mathrm{div}}\,,
\label{div1}
\end{equation}
where
\begin{equation}
\Big({\mathrm{div}}\left(a_k\right)\Big)^{i_1\dots i_{k-1}}
=
\partial_j
(a_k^{i_1\dots i_{k-1}j})\,.
\label{div}
\end{equation}

\medskip
\noindent
{\bf Remarks}.
(a)
The operator (\ref{diver})
\textit{commutes with the action of the affine Lie algebra}.
Moreover, it follows from the classical Weyl-Brauer theorem
that
any linear operator on $\cS$ commuting with the affine Lie algebra
is a polynomial expression in $\Div$ and $\rE$ defined by (\ref{rE}).

(b)
Fix the standard volume form $\Omega=dx^1\wedge\cdots\wedge dx^n$
on $\bbR^n$, then, for any vector field $X$,
the standard divergence $\mathrm{div}_{\Omega}(X)$
coincides with $\Div\,X$ (after the identification (\ref{Identific})).

%%%%%%%%%%%%%%%%%%%%%%%%%%%%%%%%%%%%%%%%%%%%%%%%%%%%%%%%%%%%%%%%%%%%%%%%%%%%%%%%%%%%
\subsection{Statement of the main result}
%%%%%%%%%%%%%%%%%%%%%%%%%%%%%%%%%%%%%%%%%%%%%%%%%%%%%%%%%%%%%%%%%%%%%%%%%%%%%%%%%%%%

The main result of this paper is the existence and uniqueness of an
$\Sl_{n+1}$-equivariant symbol map on $\bbR^n$.
We will prove that such a map is necessarily differential.

\begin{thm}
For every $\l$, there exists a unique
$\Sl_{n+1}$-equivariant symbol map $\s_\l$.
It is differential and given on a homogeneous component $\cS^k$ by
\begin{equation}
\s_\l\Big|_{\cS^k}
=
\sum_{\ell\leq k}
c^k_\ell\;{\mathrm{div}}^{k-\ell}
\label{symbol}
\end{equation}
where the numbers $c^k_\ell$ are as follows:
\begin{equation}
c_\ell^k=(-1)^{k-\ell}
\frac{{k \choose \ell}{(n+1)\lambda+k-1 \choose k-\ell}}{{k+\ell+n
\choose k-\ell}}
\label{exp}
\end{equation}
\label{theorem1}
\end{thm}
Recall that the binomial coefficient ${\a\choose m}$
for $\a\not\in\bbN$ is 
$
\frac{\a(\a-1)\cdots(\a-m+1)}{m!}.
$

\goodbreak
%\medskip
\noindent
{\bf Remarks}.

(a)
The condition of $\Sl_{n+1}$-equivariance is already sufficient to determine
$\s_\l|_{\cS^k}$ up to a constant.
The supplementary condition, that the higher order term
of the polynomial $\s_\l(A)$
coincides with the principal symbol $\s_A$, fixes the normalization.

(b) In the one-dimensional case ($n=1$)
the formula (\ref{symbol},\ref{exp}) is also valid  for the space of
pseudodifferential operators. In this case, it coincides with
(4.11) of \cite{cmz}.

%%%%%%%%%%%%%%%%%%%%%%%%%%%%%%%%%%%%%%%%%%%%%%%%%%%%%%%%%%%%%%%%%%%%%%%%%%%%%%%%%%%%
\subsection{Example: second order operators and quadratic Hamiltonians}
%%%%%%%%%%%%%%%%%%%%%%%%%%%%%%%%%%%%%%%%%%%%%%%%%%%%%%%%%%%%%%%%%%%%%%%%%%%%%%%%%%%%

Let us apply the general formul{\ae} (\ref{symbol},\ref{exp}) to
second order differential operators.

The $\Sl_{n+1}$-equivariant symbol map associates
to a second-order differential
operator
$A=a_2^{ij}\partial_i\partial_j+a_1^i\partial_i+a_0$ 
the polynomial
$
\s_\l(A)=\bar a_2^{ij}\xi_i\xi_j+\bar a_1^i\xi_i+\bar a_0,
$
with the coefficients~:
$$
\matrix{
\bar a_2^{ij}
=&
a_2 ^{ij}\hfill \cr\noalign{\smallskip}
\bar a_1^i
=& \displaystyle
a_1^i- 2\frac{(n+1)\lambda+1}{n+3}\partial_j(a_2^{ij})\hfill
\cr\noalign{\smallskip}
\bar a_0
=& \displaystyle
a_0 - \lambda \partial_i(a_1^i) +
\lambda\frac{(n+1)\lambda+1}{n+2}\partial_i\partial_j(a_2^{ij})\,.
\hfill \cr
}
$$

%%%%%%%%%%%%%%%%%%%%%%%%%%%%%%%%%%%%%%%%%%%%%%%%%%%%%%%%%%%%%%%%%%%%%%%%%%%%%%%%%%%%
\subsection{Proof of Theorem \ref{theorem1}}
%%%%%%%%%%%%%%%%%%%%%%%%%%%%%%%%%%%%%%%%%%%%%%%%%%%%%%%%%%%%%%%%%%%%%%%%%%%%%%%%%%%%

Let us first prove that symbol map $\s_\l$
defined by the formul{\ae} (\ref{symbol}) and (\ref{exp})
is $\Sl_{n+1}$-equivariant.

A symbol map (\ref{symbol}) with arbitrary constants $c^k_{\ell}$,
is obviously equivariant with respect to the affine algebra.
It can be rewritten in the form~:
\begin{equation}
\s_\l\Big|_{\cS^k}=
\sum_{\ell=O}^kC_\ell^k\,\Div^{k-\ell}\,,
\label{algebraic}
\end{equation}
where $c_\ell^k=k(k-1)\cdots(k-\ell+1)C_\ell^k$.
Now, to determine the constants $C_\ell^k$,
one needs the following commutation relations.

\begin{lem}
For $X_s=x^s\cE\in \Sl(n+1,\bbR)$, one has~:
\begin{equation}
[L_{X_s},\Div]
=
\Big(2\rE+(n+1)\Big)\circ\partial_{\xi_s}\,,
\label{relation}
\end{equation}
\label{commutation}
\end{lem}
\noindent
\textit{Proof}.
Straightforward.
\cqfd

\medskip

Let us calculate a recurrent relation for the coefficients
$C_\ell^k$ in (\ref{algebraic}).
The condition of $\Sl(n+1,\bbR^n)$-equivariance reads~:
$$
L_X\circ\s_\l=
\s_\l\circ\cL_X^\l\,,
\quad
\hbox{for every}
\quad
X\in\Sl(n+1,\bbR^n)\,.
$$
For $X=X_s$, from (\ref{SLAction})
one, therefore, gets~:
\begin{equation}
[L_{X_s},\s_\l]=
-\s_\l\circ\Big(\rE+\l(n+1)\Big)\circ\partial_{\xi_s}.
\label{TheEquation}
\end{equation}
We are looking for the symbol map $\s_\l$ in the form (\ref{algebraic}).
From the equation (\ref{algebraic},\ref{TheEquation})
one immediately obtains~:
$$
C_\ell^k\,\left[L_X,\Div^{k-\ell}\right]\Big|_{\cS^k}
=
-C_\ell^{k-1}\,
\Div^{k-\ell-1}\circ\Big(\rE+\l(n+1)\Big)\circ
\partial_{\xi_s}.
$$
From the commutation relation (\ref{relation}) and using that
$\rE=k\,\mathrm{Id}$ on $\cS^k$, one easily gets~:
$$
\left[L_{X_s},\Div^m\right]\Big|_{\cS^k}
=
m(2k-m+n)\,\Div^{m-1}\circ\partial_{\xi_s}\,,
$$
and, finally, one obtains the recurrent relation~:
\begin{equation}
C_\ell^k=
-\frac{k-1+\l(n+1)}{(k-\ell)(k+\ell+n)}C_\ell^{k-1}\,.
\label{Recurrent}
\end{equation}
This together with $C_k^k=1$ is equivalent to (\ref{exp}).

This proves that the symbol map $\s_\l$ is, indeed,
$\Sl_{n+1}$-equivariant.

\medskip

The uniqueness of the $\Sl_{n+1}$-equivariant
map $\s_\l$ follows from~:

\begin{pro}
The only $\Sl_{n+1}$-equivariant linear maps on the space of
symmetric contravariant tensor fields are multiplications by
constants, namely~:
\begin{equation}
\mathrm{Hom}_{\Sl_{n+1}}(\cS^k,\cS^\ell)
=
\left\{
\begin{array}{cc}%{r,l}
\bbR\,,&k=\ell\\
0\,,&k\not=\ell
\end{array}
\right.
\label{HomSl}
\end{equation}
\label{simple}
\end{pro}

\noindent
\textit{Proof}.
It is easy to compute the Casimir operator of $\Sl(n+1,\bbR)$~:
\begin{equation}
C\Big|_{\cS^k}
=
2k(k+n)
\,\mathrm{Id}\,.
\label{Casimir}
\end{equation}
Since $C$ commutes with the $\Sl(n+1,\bbR)$-action,
this implies the result (\ref{HomSl}) in the case $k\not=\ell$.

In the case $k=\ell$, the result easily follows from Theorem~\ref{discuss}
below. Moreover, Lemma~\ref{LEM52} implies that 
an $\Sl_{n+1}$-equivariant linear operator on $\cS^k$
is necessarily a differential operator of order zero.
Then, it follows from the classical Weyl-Brauer theorem that it
is proportional to the identity.
\cqfd

\medskip

Theorem \ref{theorem1} is proven.
\cqfd

%%%%%%%%%%%%%%%%%%%%%%%%%%%%%%%%%%%%%%%%%%%%%%%%%%%%%%%%%%%%%%%%%%%%%%%%%%%%%%%%%%%%
\subsection{Projectively equivariant quantization}
%%%%%%%%%%%%%%%%%%%%%%%%%%%%%%%%%%%%%%%%%%%%%%%%%%%%%%%%%%%%%%%%%%%%%%%%%%%%%%%%%%%%

The inverse to the symbol map, $\s_\l^{-1}$, is the unique
(up to normalization) $\Sl_{n+1}$-equivariant quantization map.

\begin{pro}
The map
$\s_\l^{-1}$ is defined by~:
\begin{equation}
\s_\l^{-1}\Big|_{\cS^k}=
\sum_{\ell\leq k}\bar c_\ell^k
\;{\mathrm{div}}^{k-\ell}
\label{TheQuantization}
\end{equation}
where
\begin{equation}
\bar c_\ell^k=
\frac{{k \choose \ell}{(n+1)\lambda+k-1 \choose k-\ell}}{{2k+n-1
\choose k-\ell}}\,.
\label{exp'}
\end{equation}
\label{theorem'}
\end{pro}

\noindent
\textit{Proof}.
Applying the same method as in the  proof of Theorem \ref{theorem1},
one can show that there exists a unique differential operator
$Q_\l:\Pol({T^*}\bbR^n)\to\cD_\l$
which is $\Sl(n+1,\bbR)$-equivariant and such that the
principal symbol of $Q_\l(P)$ coincides with the leading term of $P$,
namely, $Q_\l=\s_\l^{-1}$.
On the space of homogeneous polynomials
this map is given by
\begin{equation}
Q_\l\Big|_{\cS^k}
=
\sum_{\ell\leq k}
\bar C^k_\ell\,\Div^{k-\ell}\,.
\label{TheQuantization'}
\end{equation}
The numbers $\bar C^k_\ell$ are defined by the relation:
$$
\left\{
\matrix{
\bar C_\ell^k
& = &
\displaystyle
\frac{\ell+(n+1)\lambda}{(k-\ell)(k+\ell+n)} \;
\bar C_{\ell+1}^k
\hfill\cr\noalign{\medskip}
\bar C_k^k&=&1
\hfill\cr
}
\right.
$$
The result follows.
\cqfd

%%%%%%%%%%%%%%%%%%%%%%%%%%%%%%%%%%%%%%%%%%%%%%%%%%%%%%%%%%%%%%%%%%%%%%%%%%%%%%%%%%%%
\subsection{Example: quantization of the geodesic flow}
%%%%%%%%%%%%%%%%%%%%%%%%%%%%%%%%%%%%%%%%%%%%%%%%%%%%%%%%%%%%%%%%%%%%%%%%%%%%%%%%%%%%

Consider a nondegenerate
quadratic form $H=g^{ij}\xi_i\xi_j$
and apply the $\Sl_{n+1}$-equivariant
quantization map in the special case of $\lambda=1/2$.
It is easy to check that the result of quantization
(\ref{TheQuantization}) is a \textit{Laplace-Beltrami operator}.
Namely~:
\begin{equation}
Q_{1/2}(H)=
\Delta_g+
\frac{(n+1)}{4(n+2)}\partial_{ij}g^{ij}\,,
\label{geodesic1}
\end{equation}
where $\Delta_g=g^{ij}(\partial_i\partial_j-\Gamma_{ij}^k\partial_k)$
is the standard Laplace operator corresponding to
the metric $g=H^{-1}=g_{ij}dx^idx^j$.

Recall that in a neighborhood of each point $u$, there exist
so-called {\it normal coordinates}, that is coordinates 
characterized by:
$\Gamma_{ij}^k(u)=0$
and
$
\partial_l(\Gamma_{ij}^k)={1\over3}(R^k_{li,j}+R^k_{lj,i})\,.
$
In these coordinates,
the potential of (\ref{geodesic1}) is proportional to the scalar curvature:
$\partial_{ij}g^{ij}=(1/3)R$.

One obtains the following result. 

\begin{pro}

In the case when the normal coordinates of $g$
are compatible with the projective structure, the quantum hamiltonian
is given by
\begin{equation}
Q_{1/2}(H)=
\Delta_g+
\frac{(n+1)}{12(n+2)}R\,.
\label{geodesic}
\end{equation}
\label{quageod}
\end{pro}

\medskip
\noindent
\textbf{Remark}.
The problem of quantization of the geodesic flow on a
(pseudo)-Riemannian manifold have already been
considered by many authors (see \cite{det} end references therein).
Various methods leads to formul{\ae} of the type (\ref{geodesic})
but with different values of the multiple in front of the scalar
curvature.
The formul{\ae} (\ref{geodesic1},\ref{geodesic})
is a new version of the quantization of the geodesic flow.

%%%%%%%%%%%%%%%%%%%%%%%%%%%%%%%%%%%%%%%%%%%%%%%%%%%%%%%%%%%%%%%%%%%%%%%%%%%%%%%%%%%%
%%%%%%%%%%%%%%%%%%%%%%%%%%%%%%%%%%%%%%%%%%%%%%%%%%%%%%%%%%%%%%%%%%%%%%%%%%%%%%%%%%%%
\section{Equivariance and locality}\label{equi}
%%%%%%%%%%%%%%%%%%%%%%%%%%%%%%%%%%%%%%%%%%%%%%%%%%%%%%%%%%%%%%%%%%%%%%%%%%%%%%%%%%%%
%%%%%%%%%%%%%%%%%%%%%%%%%%%%%%%%%%%%%%%%%%%%%%%%%%%%%%%%%%%%%%%%%%%%%%%%%%%%%%%%%%%%

Consider the Lie algebra $\bbR\ltimes\bbR^n$
generated by the vector fields
\begin{equation}
\cE=x^i\frac{\partial}{\partial x^i}\,,
\quad
\frac{\partial}{\partial x^1}\,,\quad
\ldots\,,\quad
\frac{\partial}{\partial x^1}\,,
\label{EulAlg}
\end{equation}
(which is, of course, a subalgebra of $\Sl(n+1,\bbR)$).

The following result seems to be quite unexpected.

\goodbreak

\begin{thm}
If $p\geq q$, then every 
linear map $T:{\cal S}^p\to{\cal S}^q$
equivariant with respect to the Lie algebra
generated by the vector fields (\ref{EulAlg}) is local.
\label{discuss}
\end{thm}

\noindent
{\it Proof}. 
Assume that $P\in{\cal S}^p$ vanishes in a neighborhood of some
point $u\in\bbR^n$. 
We will prove that $T(P)_u=0$. Using invariance with respect
to translations $\bbR^n\subset\Sl(n+1,\bbR)$, we may assume that $u=0$.

\begin{lem}
For every $P\in{\cal S}^p$ such that the $(s-1)$-jet of $P$
at the origin vanishes, there
exists $Q\in{\cal S}^p$ such that
$$
P
=
x^{i_1}\cdots x^{i_s}\partial_{i_1}\cdots\partial_{i_s}(Q)
$$
\label{LEM52}
\end{lem}

\noindent
{\it Proof of the lemma}.
It follows directly from the Taylor integral formula.
\cqfd

\medskip

In terms of $\cE=x^i\partial_i$, this reads:
$$
P=
(\cE-(s-1)\,{\rm Id})(\cE-(s-2)\,{\rm Id})\cdots(\cE-{\rm Id})\cE (Q).
$$
Taking  account of the fact that the Lie derivative $L_\cE$ on ${\cal S}^p$ is
given by $L_\cE=\cE-p\,{\rm Id}$, we thus obtain:
$$
P=
(L_\cE+(p-s+1)\,{\rm Id})(L_\cE+(p-s+2)\,{\rm Id})\cdots(L_\cE+(p-1)\,
{\rm Id})(L_\cE+p\,{\rm Id}) (Q).
$$
Therefore, by $\Sl_{n+1}$-equivariance,
$$
\matrix{
T(P) &=&
(L_\cE+(p-s+1)\,{\rm Id})\cdots(L_\cE+p\,{\rm Id})T(Q)
\hfill\cr\noalign{\medskip}
&=&
(\cE+(p-q-s+1)\,{\rm Id})\cdots(\cE+(p-q)\,{\rm Id})T(Q).\hfill\cr
}
$$
With $s=p-q+1$, we obtain $T(P)_0=0$.
Hence the result.
\cqfd

\medskip
\noindent
{\bf Remarks}.

(a)
By the well-known Peetre theorem \cite{pee}, the map
$T$ from Theorem \ref{discuss} is (locally) a differential operator.
This theorem therefore generalizes the statement of Theorem \ref{theorem1}
that every $\Sl_{n+1}$-equivariant symbol map is
differential.

(b)
The Lie algebra generated by the vector fields (\ref{EulAlg})
is also a subalgebra of the conformal Lie algebra $\so(p+1,q+1)$,
and, therefore, can be implied in this case.

%%%%%%%%%%%%%%%%%%%%%%%%%%%%%%%%%%%%%%%%%%%%%%%%%%%%%%%%%%%%%%%%%%%%%%%%%%%%%%%%%%%%
%%%%%%%%%%%%%%%%%%%%%%%%%%%%%%%%%%%%%%%%%%%%%%%%%%%%%%%%%%%%%%%%%%%%%%%%%%%%%%%%%%%%
\section{Action of $\Vect(\bbR^n)$ on $\cD_\l$
in $\Sl_{n+1}$-equivariant form}
%%%%%%%%%%%%%%%%%%%%%%%%%%%%%%%%%%%%%%%%%%%%%%%%%%%%%%%%%%%%%%%%%%%%%%%%%%%%%%%%%%%%
%%%%%%%%%%%%%%%%%%%%%%%%%%%%%%%%%%%%%%%%%%%%%%%%%%%%%%%%%%%%%%%%%%%%%%%%%%%%%%%%%%%%

We now turn to study the space of differential operators on $\bbR^n$ as a
module over the Lie algebra $\Vect(\bbR^n)$ of all vector fields on $\bbR^n$.
Using the $\Sl_{n+1}$-equivariant symbol $\s_\l$ one
can obtain a canonical (projectively invariant) form of the
$\Vect(\bbR^n)$-action on $\cD_\l$~:
\begin{equation}
\begin{CD}
\cD_\l @>\cL^{\lambda}_X >> \cD_\l \strut\\
@V{\s_\l}VV @VV{\s_\l}V \strut\\
\Pol(T^*{}\bbR^n) @>\s_\l\circ \cL^{\lambda}_X\circ\s_\l^{-1} >>
\Pol(T^*{}\bbR^n)\strut
\end{CD}
\label{TheDiagram}
\end{equation}
This action is clearly of the form
\begin{equation}
\s_\l\circ \cL^{\lambda}_X\circ\s_\l^{-1}
=
L_X+
\sum_{\ell\geq1}
\g_\ell^\l(X)\,,
\label{nac}
\end{equation}
where $\gamma_\ell^{\lambda}$ are linear maps associating
to $X$ some operators on $\Pol(T^*{}\bbR^n)$, namely, for every $k$,
\begin{equation}
\g^\l_\ell:
\Vect(\bbR^n)\to
\Hom(\cS^k,\cS^{k-\ell})\,.
\label{gn}
\end{equation}
The operators $\g_\ell^\l$ play important role in the study of
these modules.

%%%%%%%%%%%%%%%%%%%%%%%%%%%%%%%%%%%%%%%%%%%%%%%%%%%%%%%%%%%%%%%%%%%%%%%%%%%%%%%%%%%%
\subsection{Elementary properties of the maps $\g_\ell^\l$}\label{Element}
%%%%%%%%%%%%%%%%%%%%%%%%%%%%%%%%%%%%%%%%%%%%%%%%%%%%%%%%%%%%%%%%%%%%%%%%%%%%%%%%%%%%

The following lemma is an immediate corollary of the $\Sl_{n+1}$-equivariance
of $\s_\l$~:

\begin{lem}
The linear maps (\ref{gn}) satisfy the following properties~:
\hfill\break
(a) $\Sl_{n+1}$-equivariance:
\begin{equation}
[L_X,\g_\ell^\l(Y)]
=
\g_\ell^\l([X,Y]),\quad X\in\Sl_{n+1}\,,
\label{equi}
\end{equation}
(b) vanishing on $\Sl_{n+1}$:
\begin{equation}
\g_\ell^\l(X)= 0,\quad X\in \Sl_{n+1}.
\label{nul}
\end{equation}
\label{PropGam}
\end{lem}

The next statement follows directly from the fact that (\ref{nac}) is
an action of $\Vect(\bbR^n)$~:

\begin{lem}
For every $\l$, the map $\g_1^\l$ is a 1-cocycle on $\Vect(\bbR^n)$, i.e.,
it satisfies the relation~:
\begin{equation}
[L_X,\g^\l_1(Y)]-[L_Y,\g^\l_1(X)]
=
\g^\l_1([X,Y])\,,
\label{CocEq}
\end{equation}
for every $X,Y\in\Vect(\bbR^n)$.
\label{ActLem}
\end{lem}

Moreover, in the case when $\g_1^\l$ vanishes, 
then the cocycle property holds true for $\g_2^\l$.

\medskip
\noindent
\textbf{Remark}.
These simple facts explain how
the modules of differential operators
are related to
so-called 
$\Sl(n+1,\bbR)$-relative cohomology of $\Vect(\bbR^n)$
(that is, $\Vect(\bbR^n)$-cohomology vanishing on $\Sl(n+1,\bbR)$).
This subject will be treated in a subsequent article.

\medskip

We will need explicit formul{\ae} for $\gamma_1^\lambda$ and
$\gamma_2^\lambda$.

%%%%%%%%%%%%%%%%%%%%%%%%%%%%%%%%%%%%%%%%%%%%%%%%%%%%%%%%%%%%%%%%%%%%%%%%%%%%%%%%%%%%
\subsection{Computing $\g_1^\l$}
%%%%%%%%%%%%%%%%%%%%%%%%%%%%%%%%%%%%%%%%%%%%%%%%%%%%%%%%%%%%%%%%%%%%%%%%%%%%%%%%%%%%

\begin{pro}
One has~:
\begin{equation}
\gamma^{\lambda}_1\Big|_{\cS^k}
=
\frac{(n+1)}{2(2k+n-1)}\,(2\lambda-1)\,
\ell_k\,,
\label{g1}
\end{equation}
where $\ell_k:\Vect(\bbR^n)\to\Hom(\cS^k,\cS^{k-1})$
is the following operator~:
\begin{equation}
\ell_k(X)
=
\partial_{ij}X\,\partial_{\xi_i\xi_j}-
\frac{2(k-1)}{n+1}\,
(\partial_i\circ\Div)X\,
\partial_{\xi_i}\,.
\label{ell}
\end{equation}
\label{pro33}
\end{pro}

\noindent
{\it Proof}.
Applying explicit formul{\ae}
(\ref{TheQuantization'}), (\ref{compar}) and (\ref{algebraic})
for each factor in
$\s_\l\circ\cL_X^\l\circ\s_\l^{-1}$,
one immediately obtains~:
$$
\begin{array}{ccl}
\g^\l_1(X)\Big|_{\cS^k}
&=&
\bar C_{k-1}^k\,L_X\circ\Div
+
C_{k-1}^k\,\Div\circ L_X\\[8pt]
&&
-
\frac{1}{2}\,\partial_{ij}X\,\partial_{\xi_i\xi_j}
-
\l\,(\partial_i\circ\Div)X\,\partial_{\xi_i}\\[8pt]
&=&
-(C_{k-1}^k-\frac{1}{2})\,\partial_{ij}X\,\partial_{\xi_i\xi_j}
-(C_{k-1}^k-\l)\,(\partial_i\circ\Div)X\,\partial_{\xi_i}\,.
\end{array}
$$
The formula (\ref{g1},\ref{ell}) follows.
\cqfd

\medskip
\noindent
\textbf{Important Remark}.
The operator (\ref{ell}) is, in fact, an operator of contraction
with the symmetric $(2,1)$-tensor field: 
$\ell_k(X)(P)=\langle\bar\ell(X),P\rangle$, where~:
\begin{equation}
\bar\ell(X)_{ij}^h
=
\partial_{ij}X^h
-\frac{1}{n+1}
\left(\d^h_i\partial_j+\d^h_j\partial_i\right)\partial_lX^l
\label{schwar}
\end{equation}
This expression is a 1-cocycle
on $\Vect(\bbR^n)$ vanishing on the subalgebra
$\Sl(n+1,\bbR)$ (cf.~\cite{hab,mm} and references therein).

%%%%%%%%%%%%%%%%%%%%%%%%%%%%%%%%%%%%%%%%%%%%%%%%%%%%%%%%%%%%%%%%%%%%%%%%%%%%%%%%%%%%
\subsection{Computing $\g_2^\l$}
%%%%%%%%%%%%%%%%%%%%%%%%%%%%%%%%%%%%%%%%%%%%%%%%%%%%%%%%%%%%%%%%%%%%%%%%%%%%%%%%%%%%

\begin{pro} The map $\gamma_2^\lambda$ in (\ref{nac}) is
given by
$$
\gamma^{\lambda}_2\Big|_{\cS^k}
=
\frac{1}{2(2k+n-2)(2k+n-3)}\,
s^\l_k,
$$
where the operator $s^\l_k:\Vect(\bbR^n)\to\Hom(\cS^k,\cS^{k-2})$ is~:
\begin{equation}
\begin{array}{ccl}
s^\l_k(X)
&=&\;\;
\alpha_1\,
\partial_{hij}X\,
\partial_{\xi_h\xi_i\xi_j}
+\alpha_2\,
\partial_{ij}X\,\Div\circ\partial_{\xi_i\xi_j}\\[8pt]
&&+
\beta_1\,
(\partial_{ij}\circ\Div)X\,
\partial_{\xi_i\xi_j}
+\beta_2\,
(\partial_{i}\circ\Div) X\,\Div\circ
\partial_{\xi_i}\\[8pt]
&&+
\delta\,
\partial_{ij}\partial_{\xi_h}X\,
\partial_h\partial_{\xi_i\xi_j}
\end{array}
\label{gamma2for}
\end{equation}
with the numerical coefficients~:

\begin{equation}
\matrix{
\displaystyle
\alpha_1=
-\Big((n+1)^2\lambda(\lambda-1)+\frac{1}{3}(k^2+kn+n^2-k+n)\Big)
\hfill\cr\noalign{\bigskip}
\displaystyle
\alpha_2=-\frac{2(n+1)^2\lambda(\lambda-1)+
2k^2+2kn-4k+n^2-n+2}{2k+n-1}
\hfill\cr\noalign{\bigskip}
\displaystyle\beta_1=
\frac{(4k+n-5)(n+1)\lambda(\lambda-1)-(k-2)(k-1)}{2k+n-1}
\hfill\cr\noalign{\bigskip}
\beta_2=(4k-6)(n+1)\lambda(\lambda-1)+(k-2)n
\hfill\cr\noalign{\bigskip}
\delta_2=-(n+1)^2\lambda(\lambda-1)+(k-2)(k+n-1)
\hfill\cr
}
\label{gamma2}
\end{equation}
\end{pro}

\medskip

\noindent
{\it Proof.} 
A quite complicated straightforward computation
similar to that of the proof of Proposition \ref{pro33}.
\cqfd

%%%%%%%%%%%%%%%%%%%%%%%%%%%%%%%%%%%%%%%%%%%%%%%%%%%%%%%%%%%%%%%%%%%%%%%%%%%%%%%%%%%%
\subsection{Example: case~ $k=2$}
%%%%%%%%%%%%%%%%%%%%%%%%%%%%%%%%%%%%%%%%%%%%%%%%%%%%%%%%%%%%%%%%%%%%%%%%%%%%%%%%%%%%

For $k=2$, one verifies that
the operator (\ref{gamma2for}) is as follows:
$$
s_2^\l(X)
=
-(n+1)(n+2)\,\l(\l-1)\,\bar s\,,
$$
where the operator $\bar s$ can be written in terms of the operator
(\ref{ell})~:
\begin{equation}
\bar s(X)
=
\partial_{\xi_i}\left(\ell_2(X)\right)\,\partial_i
-\frac{2}{n-1}\,\Div\left(\ell_2(X)\right)\,.
\label{g2}
\end{equation}

\medskip

\noindent
\textbf{Remark}.
In the case $\l=1/2$, the map $\g_1$
is identically zero and
the first non-zero term, namely $\g_2$, in (\ref{nac}) is a 1-cocycle
on $\Vect(\bbR^n)$ (cf. Section~\ref{Element}).

%%%%%%%%%%%%%%%%%%%%%%%%%%%%%%%%%%%%%%%%%%%%%%%%%%%%%%%%%%%%%%%%%%%%%%%%%%%%%%%%%%%%
\subsection{Conjugation and the $\Sl_{n+1}$-equivariant symbol map}
%%%%%%%%%%%%%%%%%%%%%%%%%%%%%%%%%%%%%%%%%%%%%%%%%%%%%%%%%%%%%%%%%%%%%%%%%%%%%%%%%%%%

Recall from \cite{duv} and \cite{lec} that the conjugation 
is an isomorphism of $\Vect(\bbR^n)$-modules
$$
*:{\cal D}_\l
\cong{\cal D}_{1-\l}\,,
$$
(which also exists for an arbitrary manifold $M$, see below).
It is characterized by $*({\rm Id})={\rm Id}$
and by:
\begin{equation}
*(L^\lambda_X\circ A)=-*(A)\circ L^{1-\lambda}_X
\label{star1}
\end{equation}
for all $X\in\Vect(\bbR^n)$ and all $A\in{\cal D}_\l$.

A nice property of the $\Sl_{n+1}$-equivariant
symbol map is the following

\begin{lem}
 For each homogeneous polynomial $P\in{\cal S}^k$ one has~:
$
\s_\l\circ *\circ\s_\l^{-1}(P)
=
(-1)^kP\,.
$
\label{lem34}
\end{lem}

\noindent
{\it Proof}. By Proposition \ref{simple}, the map
$\s_\l\circ *\circ\s_\l^{-1}$ is an operator of
multiplication by a constant.
On the other hand, (\ref{star1}) shows that the principal symbol
$\s_{*(A)}=(-1)^k\s_A$.
\cqfd

\begin{cor} 
All the bilinear maps
$\gamma_{2p+1}^{\frac{1}{2}}$ vanish.
\label{pro35}
\end{cor}

%%%%%%%%%%%%%%%%%%%%%%%%%%%%%%%%%%%%%%%%%%%%%%%%%%%%%%%%%%%%%%%%%%%%%%%%%%%%%%%%%%%%
%%%%%%%%%%%%%%%%%%%%%%%%%%%%%%%%%%%%%%%%%%%%%%%%%%%%%%%%%%%%%%%%%%%%%%%%%%%%%%%%%%%%
\section{Space of differential operators on a manifold
as a module over the Lie algebra of vector fields}
%%%%%%%%%%%%%%%%%%%%%%%%%%%%%%%%%%%%%%%%%%%%%%%%%%%%%%%%%%%%%%%%%%%%%%%%%%%%%%%%%%%%
%%%%%%%%%%%%%%%%%%%%%%%%%%%%%%%%%%%%%%%%%%%%%%%%%%%%%%%%%%%%%%%%%%%%%%%%%%%%%%%%%%%%

Let $M$ be a smooth manifold.
Consider the space $\cF_\l$ of tensor densities of degree $\l$ on $M$ and
${\cal D}_\l(M)$ be the space of scalar linear differential operators
$A:\cF_\l\to\cF_\l$.
The space ${\cal D}_\l(M)$ has a natural $\Vect(M)$-module structure.
The 1-parameter family of $\Vect(M)$-modules $\cD_\l$ has been recently
studied in \cite{duv},\cite{lec},\cite{gar}.

In this section we will study the quotient modules
${\cal D}^k_{\lambda}(M)/{\cal D}^\ell_{\lambda}(M)$.
Our purpose is to solve the problem of isomorphism between these modules
for different values of $\lambda$
and to compare these modules with the modules
of symmetric contravariant tensor fields
\allowbreak
$
\Pol^k(T^*M)/\Pol^\ell(T^*M)
=
{\cal S}^k(M)\oplus\dots\oplus{\cal S}^{\ell+1}(M)$.

%%%%%%%%%%%%%%%%%%%%%%%%%%%%%%%%%%%%%%%%%%%%%%%%%%%%%%%%%%%%%%%%%%%%%%%%%%%%%%%%%%%%
\subsection{Locality of equivariant maps}
%%%%%%%%%%%%%%%%%%%%%%%%%%%%%%%%%%%%%%%%%%%%%%%%%%%%%%%%%%%%%%%%%%%%%%%%%%%%%%%%%%%%

Equivariance with respect to $\Vect(M)$ is of course a much stronger
condition that equivariance with respect to $\Sl(n+1,\bbR)$
or $\bbR\ltimes\bbR^n$ (cf. Section \ref{equi}).
It has already been shown in \cite{lec} that a $\Vect(M)$-equivariant
linear map from ${\cal D}^k_{\lambda}(M)$ into ${\cal D}^k_{\mu}(M)$
is local. In the same spirit we have:

\begin{pro}
Every $\Vect(M)$-equivariant
linear map $T$ from ${\cal D}^k_{\lambda}(M)/{\cal D}^\ell_{\lambda}(M)$
into ${\cal D}^k_{\mu}(M)/{\cal D}^\ell_{\mu}(M)$ is local ($k>\ell$).
\label{pro41}
\end{pro}

{\bf Proof.}
Suppose that $A\in{\cal D}^k_{\lambda}(M)/{\cal D}^\ell_{\lambda}(M)$
vanishes on an open subset $U$ of $M$. Let us show that $T(A)$
vanishes on $U$ as well. Assume in the contrary that $T(A)|_u\not=0$ for
some $u\in U$. One can choose $u$ such that the principal
symbol $\s_{T(A)}|{}_u\not=0$.
Then, there exists $X\in\Vect(M)$ with compact support in $U$
such that $(L_X(T(A)))|_u\not=0$. Hence the contradiction since $L_X(A)=0$
and $L_X(T(A))=T(L_X(A))$. To see that such $X$ exists, it suffice
to choose $X$ such that $(L_X(S))|_u\not=0$ where
$S\in{\cal S}^p(M),\,\ell<p\leq k$, is the principal symbol of $T(A)$
and $p$ is its order around $u$. This is always possible for $p>0$.
Hence the proposition.
\cqfd

\medskip

The above proof is easily adapted to get the following result that
we need for later purpose.

\begin{pro}
 Every $\Vect(M)$-equivariant linear map
$T:{\cal D}^k_{\lambda}(M)/{\cal D}^\ell_{\lambda}(M)\to{\cal S}^p(M)$
is local ($k>\ell,p>0$).
\label{pro42}
\end{pro}

%%%%%%%%%%%%%%%%%%%%%%%%%%%%%%%%%%%%%%%%%%%%%%%%%%%%%%%%%%%%%%%%%%%%%%%%%%%%%%%%%%%%
\subsection{Classification of the
modules ${\cal D}^k_{\lambda}(M)/{\cal D}^{k-2}_{\lambda}(M)$}
%%%%%%%%%%%%%%%%%%%%%%%%%%%%%%%%%%%%%%%%%%%%%%%%%%%%%%%%%%%%%%%%%%%%%%%%%%%%%%%%%%%%

The following result gives a classification of the quotient modules
for the case $\ell=k-2$.

\begin{thm} Assume that $k \geq 2$ and that $\dim M \geq 2$.
\hfill\break
(i)All the $\Vect(M)$-modules
$
{\cal D}^k_{\lambda}(M)/{\cal D}^{k-2}_{\lambda}(M)
$,
with $\lambda\not=1/2$,
are isomorphic. They are  not isomorphic to the direct sum
of modules of tensor fields ${\cal S}^k(M) \oplus {\cal S}^{k-1}(M)$.
\hfill\break
(ii)
The module of differential operators on half-densities is exceptional since
$$
{\cal D}^k_{1\over2}(M)/{\cal D}^{k-2}_{1\over2}(M)
\cong
{\cal S}^k(M)\oplus{\cal S}^{k-1}(M).
$$
\label{thm43}
\end{thm}

\noindent
{\it Proof}.
(i) Having local coordinates in a open subset $U\subset M$ and
using the $\Sl_{n+1}$-equivariant symbol map,
we identify ${\cal D}^k_{\lambda}(M)/{\cal D}^{k-2}_{\lambda}(M)$
over $U$ with ${\cal S}^k(M)\oplus{\cal S}^{k-1}(M)$
endowed with the $\Vect(M)$-action
\begin{equation}
\cL_X^\l\Big|_{\cS^k\oplus\cS^{k-1}}
=
L_X+\gamma_1^\lambda(X)
\label{secact}
\end{equation}
In view of (\ref{g1}), it is clear that if $\lambda,\mu\not=1/2$,
the map
\begin{equation}
(P_k,P_{k-1})\mapsto
\left(P_k,\;\frac{2\mu-1}{2\lambda-1}\,P_{k-1}\right)
\label{ism}
\end{equation}
is an isomorphism between the $\Vect(\bbR^n)$-modules~:
$$
({\cal S}^k\oplus{\cal S}^{k-1},\;
\sigma_\lambda\circ\cL_X^\lambda\circ\sigma_\lambda^{-1})
\longrightarrow
({\cal S}^k\oplus{\cal S}^{k-1},\;
\sigma_\mu\circ\cL_X^\mu\circ\sigma_\mu^{-1})\,.
$$
It, therefore, defines an isomorphism between the restrictions of the
$\Vect(M)$-modules ${\cal D}^k_{\lambda}(M)/{\cal D}^{k-2}_{\lambda}(M)$
and ${\cal D}^k_{\mu}(M)/{\cal D}^{k-2}_{\mu}(M)$ to the domain $U$
of coordinates. 
This isomorphism does not depend on the choice of coordinates
because it commutes with the $\Vect(M)$-action: the formula
(\ref{ism}) does not change under the coordinate transformations.
Hence, the local isomorphisms defined on each $U$ glue together
to define a global isomorphism.

Let now ${\cal D}^k_{\lambda}(M)/{\cal D}^{k-2}_{\lambda}(M)
\to
{\cal S}^k(M)\oplus{\cal S}^{k-1}(M)$ be a $\Vect(M)$-equivariant map.
It is local (Proposition \ref{pro42}). Hence, by the Peetre theorem \cite{pee},
it is locally a differential operator. Expressed like above in terms of
the $\Sl_{n+1}$-equivariant symbol, it has a diagonal form
(Proposition \ref{simple}):
$(\bar a_k,\bar a_{k-1})\mapsto(\alpha_k\bar a_k,\alpha_{k-1}\bar a_{k-1})$
for some $\alpha_k,\alpha_{k-1}\in\bf R$. The fact that $T$ intertwines
the action (\ref{secact}) and the Lie derivative of tensors on
${\cal S}^k(M)\oplus{\cal S}^{k-1}(M)$ implies
$\alpha_{k-1}\gamma_1^\lambda(X,\bar a_k)=0$ for all
$\bar a_k\in{\cal S}^k(M)$. Since $\lambda\not=1/2$, it follows that
$\alpha_{k-1}=0$. Therefore, $T$ is not injective.

Hence Part (i) of Theorem \ref{thm43}.

%%%%%%%%%%%%%%%%%%%%%%%%%%%%%%%%%%%%%%%%%%%%%%%%%%%%%%%%%%%%%%%%%%%%%%%%%%%%%%%%%%%%
\subsection{Exceptional case $\lambda=1/2$}
%%%%%%%%%%%%%%%%%%%%%%%%%%%%%%%%%%%%%%%%%%%%%%%%%%%%%%%%%%%%%%%%%%%%%%%%%%%%%%%%%%%%

If  $\lambda=1/2$, the term
$\gamma_1^\lambda(X,\bar a_k)$ vanishes.
The $\Vect(M)$-action (\ref{nac})
in this case is just the standard action on
${\cal S}^k(M)\oplus{\cal S}^{k-1}(M)$.

Hence Theorem \ref{thm43}.
\cqfd

\medskip

It follows from this theorem that for $\lambda\not=1/2$,
there is no intrinsically defined subsymbol of
a differential operator.
However, in the exceptional case of differential operators on
$1/2$-densities, the two first terms of the
$\Sl_{n+1}$-equivariant symbol have geometric meaning. One obtains
the following remark.

\begin{cor}
 The
$\Sl_{n+1}$-equivariant symbol
defines a $\Vect(M)$-equivariant map
$$
(\sigma^{1/2}_k,\sigma^{1/2}_{k-1}):
{\cal D}^k_{1/2}(M)\to
{\cal S}^k(M)\oplus{\cal S}^{k-1}(M).
$$
\label{cor44}
\end{cor}

\noindent
\textit{Proof}.
The $\Vect(M)$-module ${\cal D}^k_{1/2}(M)$ has a symmetry:
the conjugation of operators $A \to A^\ast$ (cf. \cite{duv,lec}).
Every  $A\in{\cal D}^k_{1/2}(M)$ has a decomposition $A=A_0+A_1$, where
$A_0^*=(-1)^kA_0$ and $A_1^*=(-1)^{k-1}A_1$.
One easily sees that $\s_{1/2}{}(A)_k$ is the principal symbol of $A_0$ and
$\s_{1/2}{}(A)_{k-1}$ that of $A_1$.
\cqfd

%%%%%%%%%%%%%%%%%%%%%%%%%%%%%%%%%%%%%%%%%%%%%%%%%%%%%%%%%%%%%%%%%%%%%%%%%%%%%%%%%%%%
\subsection{Modules ${\cal D}^k_{\lambda}(M)/{\cal D}^\ell_{\lambda}(M)$
in multi-dimensional case}
%%%%%%%%%%%%%%%%%%%%%%%%%%%%%%%%%%%%%%%%%%%%%%%%%%%%%%%%%%%%%%%%%%%%%%%%%%%%%%%%%%%%

We now turn to the $\Vect(M)$-modules
${\cal D}^k_{\lambda}(M)/{\cal D}^\ell_{\lambda}(M)$
when $k-\ell\geq 3$ and $\dim M\geq 2$.
The following result shows that there is no nontrivial isomorphism
between the $\Vect(M)$-modules in this case.

\begin{thm}
 Assume that $k-\ell \geq 3$ and $\dim M \geq 2$.
\hfill\break
(i) If $\lambda\not=\mu$, the modules
${\cal D}^k_{\lambda}(M)/{\cal D}^\ell_{\lambda}(M)$
and
${\cal D}^k_{\mu}(M)/{\cal D}^\ell_{\mu}(M)$
are isomorphic if and only if $\lambda+\mu=1$.
\hfill\break
(ii) There is no isomorphism
between the modules ${\cal D}^k_{\lambda}(M)/{\cal D}^\ell_{\lambda}(M)$
and the module of tensor fields
${\cal S}^k(M)\oplus\dots\oplus{\cal S}^{\ell+1}(M)$.
\label{thm45}
\end{thm}

\noindent
{\it Proof}.
The isomorphism in Part (i) is given by the standard conjugation
of differential operators.
This map defines a general isomorphism
$
*:{\cal D}_{\lambda}(M)\to{\cal D}_{1-\lambda}(M)
$
(cf. \cite{duv,lec}).

The proof that there is no other isomorphism is very similar to that of
Theorem \ref{thm43}.
In local coordinates and
in terms of $\Sl_{n+1}$-equivariant symbol, such an isomorphism
between ${\cal D}^k_{\lambda}(M)/{\cal D}^{\ell}_{\lambda}(M)$
and ${\cal D}^k_{\mu}(M)/{\cal D}^{\ell}_{\mu}(M)$
(or $S^k(M) \oplus \cdots \oplus S^{\ell +1}(M)$)
is of the form:
\begin{equation}
(P_k,\;P_{k-1},\;P_{k-2}, \ldots)
\rightarrow
(\alpha_k\,P_k,
\;\alpha_{k-1}P_{k-1},
\;\alpha_{k-2}\,P_{k-2},\ldots )
\label{net}
\end{equation}
for some $\alpha_k,\alpha_{k-1},\alpha_{k-2},\ldots\in\bf R$.

In the both cases, (i) and (ii), the $\Vect(M)$-equivariance condition now
involves not only $\gamma_1^\lambda$ but also $\gamma_2^\lambda$.
Tedious computations allow then to show in the first case that
(\ref{net}) is an isomorphism if and only if
$(\lambda-\mu)(\lambda+\mu-1)=0$. In the second case, equivariance
immediately implies that $\alpha_{k-1}=\alpha_{k-2}=0$.
\cqfd

\medskip

Part (ii) of Theorem \ref{thm45} confirms the fact (well-known ``in practice'')
that a ``complete" symbol of a differential operator can not be defined in
an intrinsic way (even in the case of differential operators on
$1/2$-densities).

%%%%%%%%%%%%%%%%%%%%%%%%%%%%%%%%%%%%%%%%%%%%%%%%%%%%%%%%%%%%%%%%%%%%%%%%%%%%%%%%%%%%
\subsection{Modules of second order differential operators}
%%%%%%%%%%%%%%%%%%%%%%%%%%%%%%%%%%%%%%%%%%%%%%%%%%%%%%%%%%%%%%%%%%%%%%%%%%%%%%%%%%%%

Consider the modules of second order differential operators
${\cal D}^2_{\lambda}(M)$.

These modules have been classified in \cite{duv} where it has been shown
that there are
exactly three isomorphism classes of $\Vect(M)$-modules among them, namely
$$
\{\cD_\lambda^2(M),\;
\lambda \neq 0, 1/2,1\},
\quad
\{\cD_0^2(M),\cD_1^2(M)\}
\quad\hbox{and}\quad
\{\cD_{1/2}^2(M)\}\,.
$$

For $\lambda,\mu\not=0,{1\over2},1$, there exists a unique (up to a constant)
intertwining operator
$$
{\cal L}^2_{\lambda,\mu}:{\cal D}^2_{\lambda}(M)\to{\cal D}^2_{\mu}(M)
$$
This has been proven in \cite{duv} within the class of local maps. It has
been shown in \cite{lec} that each $\Vect(M)$-equivariant map from
${\cal D}^2_{\lambda}(M)\to{\cal D}^2_{\mu}(M)$ is local.

Let us express ${\cal L}^2_{\lambda,\mu}$
in terms of the $\Sl_{n+1}$-equivariant symbol (\ref{exp}).
It follows from Corollary 2.6, that
the map
$\sigma_{\mu}\circ{\cal L}^2_{\lambda,\mu}\circ\sigma_{\lambda}^{-1}$
is diagonal.

\begin{pro} The
$\Vect(M)$-module isomorphism
$\sigma_{\mu}\circ{\cal L}^2_{\lambda,\mu}\circ\sigma_{\lambda}^{-1}$
is of the form~:
\begin{equation}
(P_2,\;
P_1,\;
P_0)
\mapsto
\left(P_2,\;
\frac{2\mu-1}{2\lambda-1}\,P_1,\;
\frac{\mu(\mu-1)}{\lambda(\lambda-1)}\,P_0
\right)\,.
\label{triism}
\end{equation}
\label{pro46}
\end{pro}

\noindent
{\it Proof}. This formula follows from  the explicit expression for the
$\Vect(M)$-action (\ref{nac})
in the case of second order operators. Straightforward computations
(cf. Section 5.5) give
$$
\sigma_{\mu}\circ{\cal L}^2_{\lambda,\mu}\circ\sigma_{\lambda}^{-1}
=
L_X
+(2\lambda-1)\widetilde\gamma_1(X)
+
\lambda(\lambda-1)\widetilde\gamma_2(X)
$$
where the maps $\widetilde\gamma_1$ and $\widetilde\gamma_2$
do not depend on $\l$;
$\widetilde\gamma_2$ is defined by (\ref{g2})
according to Proposition \ref{pro33}, and
$
\widetilde\gamma_1
$
is defined by (\ref{g1},\ref{ell}).
As in the proof of Theorem \ref{thm43}, one notes that the map defined by
(\ref{triism}) is intrinsic.
\cqfd

\medskip

Let us now give a coordinate free expression for the map 
${\cal L}^2_{\lambda,\mu}$.

Every second order differential operator can be written
 as a sum of~:
\hfill\break
(i) a zero-order operator $\phi\mapsto f\phi$ (multiplication by a function),
\hfill\break
(ii) a first order operator $L^{\lambda}_X$ (Lie derivative),
\hfill\break
(iii) a symmetric expression
$
[L^{\lambda}_X,L^{\lambda}_Y]_+=
L^{\lambda}_X\circ L^{\lambda}_Y+L^{\lambda}_Y\circ L^{\lambda}_X,
$
\hfill\break
where  $f\in C^{\infty}(M),X,Y,Z\in\Vect(M)$.

\begin{pro}
 One has
$$
\begin{array}{rcl}
{\cal L}^2_{\lambda,\mu}
\Big([L^{\lambda}_X,L^{\lambda}_Y]_+\Big)&=&
[L^{\mu}_X,L^{\mu}_Y]_+ \hfill\\\noalign{\smallskip}
{\cal L}^2_{\lambda,\mu}
(L^{\lambda}_Z)&=&
\displaystyle\frac{2\mu-1}{2\lambda-1}L^{\mu}_Z\hfill\\\noalign{\smallskip}
{\cal L}^2_{\lambda,\mu}(f)&=&
\displaystyle\frac{\mu(\mu-1)}{\lambda(\lambda-1)}f\hfill\\
\end{array}
$$
\label{pro47}
\end{pro}

\noindent
{\it Proof}. Straightforward
\cqfd

\medskip
\noindent
\textbf{Remarks}.

(a)
The expression for ${\cal L}^2_{\lambda,\mu}$
in terms of Lie derivatives is intrinsic, but it is a
nontrivial fact
that it does not depend on the choice of $X,Y$ and $f$
representing the same differential operator.
The expression for ${\cal L}^2_{\lambda,\mu}$
in terms of symbols is well-defined locally, but it is
a nontrivial fact that it is invariant with respect to coordinate changes.
The two facts are corollaries of the third one: the
two formul{\ae} represent the same map.

(b)
The explicit formula for ${\cal L}^2_{\lambda,\mu}$
in terms of coefficients of differential operators
was obtained in \cite{duv}.

%%%%%%%%%%%%%%%%%%%%%%%%%%%%%%%%%%%%%%%%%%%%%%%%%%%%%%%%%%%%%%%%%%%%%%%%%%%%%%%%%%%%
\subsection{Modules of (pseudo)differential operators
in the one-dimensional case}
%%%%%%%%%%%%%%%%%%%%%%%%%%%%%%%%%%%%%%%%%%%%%%%%%%%%%%%%%%%%%%%%%%%%%%%%%%%%%%%%%%%%

We now study  the space of (pseudo)differential operators
on $S^1$ (or on $\bbR$)~:
\begin{equation}
A=\sum_{i=0}^{\infty}a_{k-i}\bigg(\frac{d}{dx}\bigg)^{k-i}\,,
\label{pse}
\end{equation}
where $k\in\bbR$.

The group $\Diff(S^1)$ and the Lie algebra $\Vect(S^1)$
act on the space of pseudodifferential operators
in the same way as on the space of differential operators.
Denote $\Psi{\cal D}^k_{\lambda}$ the
$\Vect(S^1)$-modules of the operators (\ref{pse})
acting on $\cF_\l$.

\begin{defi}
The bilinear operations
$J_m:\cF_\l\otimes\cF_\m\to\cF_{\l+\m+m}$ given by~:
\begin{equation} \label{E20}
J_m(\phi,\psi) = \sum_{i+j=m}(-1)^i m! 
%\left( \begin{array}{c}
{2\lambda+m-1 \choose
i}
{2\mu+m-1 \choose
j}
\phi^{(i)} \psi^{(j)}
\end{equation}
is called the transvectants.
\label{defitra}
\end{defi}
It is a classical fact that the transvectants are
$Sl_{n+1}$-equivariant. Moreover, for generic $\l,\m$
the operations (\ref{E20}) are characterized by this property.

Put $\l=-1$ and $\phi=X$ a vector field on $S^1$, the transvectants
(\ref{E20}) define the linear maps
\begin{equation}
\bar\g_m:\Vect(S^1)\to\Hom(\cF_\m,\cF_{\m+m-1})
\label{tra}
\end{equation}
such that $\bar\g_m(X,\psi):=J_m(X,\psi)$.

\begin{lem}
The $\Vect(S^1)$-action on $\Psi{\cal D}^k_{\lambda}$ is of the form
\begin{equation}
\s_\l\circ\cL_X\circ\s_\l^{-1}
= 
L_X
+
\sum_{i=s+2}^{k} t_{i}^{i-s}(\l)\,\bar\g_{i-s+1}(X)
\label{ber}
\end{equation}
where $t_i^{i-s}(\lambda)$ are some polynomials and $\bar\g_m$
are the operations (\ref{tra}).
\label{actlem1}
\end{lem}

\noindent
\textit{Proof}.
The formula (\ref{ber}) is a particular case of (\ref{nac}).
If $n=1$, then $\g_1\equiv0$ (cf. Proposition~\ref{pro33})
and $\g_m$ are proportional to $\bar\g_m$
since the transvectants are unique $\Sl_2$-equivariant bilinear
maps.
\cqfd

We will study the quotient-modules
$\Psi{\cal D}^k_{\lambda}/\Psi {\cal D}^{k-\ell}_{\lambda}$.
If $k\in\bbN, \ell=k+1$, it is
just the module of differential operators on $S^1$.

%%%%%%%%%%%%%%%%%%%%%%%%%%%%%%%%%%%%%%%%%%%%%%%%%%%%%%%%%%%%%%%%%%%%%%%%%%%%%%%%%
\subsection{Classification of the modules
$\Psi{\cal D}^k_{\lambda}/\Psi {\cal D}^{k-\ell}_{\lambda}$}
%%%%%%%%%%%%%%%%%%%%%%%%%%%%%%%%%%%%%%%%%%%%%%%%%%%%%%%%%%%%%%%%%%%%%%%%%%%%%%%%%

The classification of $\Vect(S^1)$-modules
$\Psi{\cal D}^k_{\lambda}/\Psi {\cal D}^{k-\ell}_{\lambda}$
follows from the formula (\ref{ber}).
As in the multi-dimensional case,
zeroes of the polynomials $t_{k-\ell}^j(\lambda)$
corresponds to exceptional modules.

Let us formulate the result
for generic values of $k$.

\begin{pro}
If $k\not=0,1/2,1,3/2,\dots$, then
the $\Vect(S^1)$-modules
$
\Psi{\cal D}^k_{\lambda}/\Psi {\cal D}^{k-\ell}_{\lambda}$
and
$
\Psi{\cal D}^k_{\mu}/\Psi {\cal D}^{k-\ell}_{\mu}
$,
are isomorphic in the following cases:
\hfill\break
(i) $\ell\geq2$;
\hfill\break
(ii) $\ell=3$, if $t^2_k(\lambda),t^2_k(\mu)\not=0$;
\hfill\break
(iii) $\ell=4$, if
$\lambda,\mu$ are not  root of the polynomials
$t^2_k,t^2_{k-1},t^3_k$;
\hfill\break
(iii) $\ell\geq4$, if and only if $\lambda+\mu=1$.
\label{pro52}
\end{pro}

The proof is analogous to the proofs of Theorems \ref{thm43} and \ref{thm45}.
\cqfd

%%%%%%%%%%%%%%%%%%%%%%%%%%%%%%%%%%%%%%%%%%%%%%%%%%%%%%%%%%%%%%%%%%%%%%%%%%%%%%%%%
\subsection{Relation to the Bernoulli polynomials}
%%%%%%%%%%%%%%%%%%%%%%%%%%%%%%%%%%%%%%%%%%%%%%%%%%%%%%%%%%%%%%%%%%%%%%%%%%%%%%%%%

Let us give the explicit formul{\ae} for the polynomials
$t^2_k,t^2_{k-1}$ and $t^3_k$~:
$$
\matrix{
\displaystyle t_k^2(\lambda)=
\frac{k(k-1)}{2k-1}\bigg(\lambda^2-\lambda-\frac{(k+1)(k-2)}{12}\bigg)
\hfill \cr\noalign{\smallskip}
\displaystyle t_k^3(\lambda)=
\frac{k}{6}\lambda(2\lambda-1)(\lambda-1) \hfill \cr
}
$$

They evoke a possible relationship between the  polynomials $t_k^j(\lambda)$
and the well-known {\it Bernoulli polynomials}. Indeed,
$$
t_k^2(\lambda)=
\frac{k(k-1)}{2k-1}\bigg(B_2(\lambda)-\frac{k(k-1)}{12}\bigg),\;\;\;
t_k^3(\lambda)=\frac{k}{12}B_3(\lambda)
$$
where $B_s$ is the Bernoulli polynomial of degree $s$, e.g.:
$$
\matrix{
B_0(x)=1,\;\;\;
B_1(x)=x-1/2,\;\;\;
B_2(x)=x^2-x+1/6,\hfill\cr
B_3(x)=x^3-3x^2/2+x/2,\;\;\;B_4(x)=x^4-2x^3+x^2-1/30,\hfill\cr
B_5(x)=x^5-5x^4/2+5x^3/3-x/6.\hfill\cr
}
$$
The next examples are:
$$
\matrix{
\displaystyle t_k^4(\lambda)=
\frac{k(k-1)(k-2)}{2(2k-3)(2k-5)}
\bigg(B_4(\lambda)+\frac{2k^2-6k+3}{24}B_2(\lambda) \hfill
\cr\noalign{\smallskip}
\displaystyle \;\;\;\;\;\;\;\;\;\;\;\;\;\;\;\;\;\;\;\;
\;\;\;\;\;\;\;\;\;\;\;\;\;\;\;\;\;\;\;\;\;\;\;\;
-\frac{3k^4+18k^3-35k^2+8k+2}{480}\bigg)\,, \hfill \cr\noalign{\smallskip}
\displaystyle t_k^5(\lambda)=
\frac{k(k-1)}{15(2k-7)}\bigg(B_5(\lambda)
+\frac{5(k-1)(k-3)}{24}B_3(\lambda)\bigg)
\hfill \cr
}
$$

\begin{pro}
Polynomials $t_k^{2j}(\lambda)$ are combinations
of $B_{2s}$ with $s=0,1,\dots,j$ and polynomials
$t_k^{2j+1}(\lambda)$ are combinations
of $B_{2s+1}$ with $s=0,1,\dots,j$
\end{pro}

\noindent
{\it Proof}. This statement is a simple corollary of the isomorphism
${\cal D}_{\lambda}\cong{\cal D}_{1-\lambda}$. Indeed, under the
involution $\lambda'=1/2-\lambda$, one has $B_s(\lambda')=(-1)^sB_s(\lambda)$.
\cqfd

\medskip
\noindent
{\bf Remark: the duality}.
There exists a nondegenerate natural pairing of
 $\Psi {\cal D}^k/\Psi {\cal D}^\ell$ and
$\Psi {\cal D}^{-\ell-2}/\Psi {\cal D}^{-k-2}$.
It is given by the so-called {\it Adler trace} \cite{adl}: if $A\in \Psi {\cal
D}^k$, where $k\in{\bf Z}$, then
$$
\hbox{tr}(A)=\int_{S^1}a_1(x)dx.
$$

Let now $A\in\Psi {\cal D}^k/\Psi \cD^\ell$ and
$B\in\Psi {\cal D}^{-\ell-2}/\Psi {\cal D}^{-k-2}$.
Put
$
(A,B):=
tr(\widetilde{A}\widetilde{B}),
$
 where
$\widetilde{A}\in\Psi {\cal D}^k,\widetilde{B}\in\Psi {\cal D}^{-\ell-2}$
are arbitrary lifts of $A$ and $B$.

Adler's trace is equivariant with respect to the $\Vect(S^1)$-action.
This means that the pairing $(\;,\;)$ is well-defined on $\Vect(S^1)$-modules.
Indeed, $([L^{\lambda}_X,A],B)+(A,[L^{\lambda}_X,B])=0$ for every
$X\in \Vect(S^1)$ (see \cite{cmz} for the details and interesting properties
of the transvectants).

%%%%%%%%%%%%%%%%%%%%%%%%%%%%%%%%%%%%%%%%%%%%%%%%%%%%%%%%%%%%%%%%%%%%%%%%%%%%%%%%%%%%
%%%%%%%%%%%%%%%%%%%%%%%%%%%%%%%%%%%%%%%%%%%%%%%%%%%%%%%%%%%%%%%%%%%%%%%%%%%%%%%%%%%%
\section{Conclusion}
%%%%%%%%%%%%%%%%%%%%%%%%%%%%%%%%%%%%%%%%%%%%%%%%%%%%%%%%%%%%%%%%%%%%%%%%%%%%%%%%%%%%
%%%%%%%%%%%%%%%%%%%%%%%%%%%%%%%%%%%%%%%%%%%%%%%%%%%%%%%%%%%%%%%%%%%%%%%%%%%%%%%%%%%%

%%%%%%%%%%%%%%%%%%%%%%%%%%%%%%%%%%%%%%%%%%%%%%%%%%%%%%%%%%%%%%%%%%%%%%%%%%%%%%%%%%%%
\subsection{Generalization for the locally projective manifolds}
%%%%%%%%%%%%%%%%%%%%%%%%%%%%%%%%%%%%%%%%%%%%%%%%%%%%%%%%%%%%%%%%%%%%%%%%%%%%%%%%%%%%

A projective structure on a manifold $M$ is defined by
an atlas
with linear-fractional coordinate changes.

More precisely, a covering $(U_i)$ with a family of local diffeomorphisms
$\phi_i:U_i\to\bbRP^n$ is called a projective atlas if the local
transformations
$\phi_j\circ\phi_i^{-1}:\bbRP^n\to\bbRP^n$ are projective (i.e. are given
by the action of the group $\SL(n+1,\bbR)$ on $\bbRP^n$).

A projective structure defines locally on $M$ an action of
the Lie group  $\SL(n+1,\bbR)$ by {\it linear-fractional transformations} and
a (locally defined) action of the Lie algebra \allowbreak $\Sl(n+1,\bbR)$
generated by vector fields (\ref{E1}),
for every system of local coordinates of a projective atlas.
This action is stable with respect to linear-fractional transformations
(the space of vector fields (\ref{E1}) is well-defined
globally on $\bbRP^n$).

 One has the following simple corollary of
Theorem \ref{theorem1}.

\begin{cor}
Given a manifold $M$ endowed with a projective
structure,
the symbol map $\s_\l$, given in local coordinates of an
arbitrary projective atlas
by the formul{\ae} (\ref{symbol},\ref{exp}), is well defined globally on $M$.
\end{cor}

%%%%%%%%%%%%%%%%%%%%%%%%%%%%%%%%%%%%%%%%%%%%%%%%%%%%%%%%%%%%%%%%%%%%%%%%%%%%%%%%%%%%
\subsection{$\Sl(n+1,\bbR)$-equivariant star-products on $T^*M$}
%%%%%%%%%%%%%%%%%%%%%%%%%%%%%%%%%%%%%%%%%%%%%%%%%%%%%%%%%%%%%%%%%%%%%%%%%%%%%%%%%%%%

Let us show that for every $\l$
the $\Sl(n+1,\bbR)$-equivariant quantization map $\s_\l$
defines a star-product on $T^*M$ thus obtaining a 1-parameter family
of $\Sl(n+1,\bbR)$-equivariant star-products.

Given a quantization map
$\sigma^{-1}:\Pol(T^*M)\to\cD(M)$,
let us introduce a new parameter $\hbar$.
For a homogeneous polynomial
$P$ of degree $k$ we set
$$
Q_{\hbar}(P)=\hbar^k\sigma^{-1}(P).
$$
and we define a new associative but non-commutative
multiplication on $\Pol(T^*M)$ by
\begin{equation}
F\star_{\hbar} G:=Q_{\hbar}^{-1}(Q_{\hbar}(F)\cdot Q_{\hbar}(G)).
\label{star}
\end{equation}
The corresponding algebra is isomorphic to the associative algebra
of differential operators on $M$.

The result of the operation (\ref{star}) is
a formal series in $\hbar$.
It has the following form:
$$
F\star_{\hbar}G=
FG+\sum_{k\geq1}\hbar^kC_k(F,G),
$$
where the higher order terms $C_k(F,G)$ are some differential operators.

\medskip

Recall that such an operation  is called a star-product
if the skew symmetric part of $C_1(F,G)$ is  the standard Poisson bracket
on ${\rm Pol}(T^*M)$.

An elementary calculation shows that the associative operation
corresponding to the quantization map (\ref{exp'}) has this property.

%%%%%%%%%%%%%%%%%%%%%%%%%%%%%%%%%%%%%%%%%%%%%%%%%%%%%%%%%%%%%%%%%%%%%%%%%%%%%%%%%
\subsection{Modules of differential operators and cohomology of $\Vect(M)$}
%%%%%%%%%%%%%%%%%%%%%%%%%%%%%%%%%%%%%%%%%%%%%%%%%%%%%%%%%%%%%%%%%%%%%%%%%%%%%%%%%

According to the formula (\ref{compar}),
the $\Vect(M)$-module of differential operators $\cD_\l(M)$, for every $\l$,
can be naturally viewed as a deformation of the module $\Pol(T^*{}\bbR^n)$ of
symmetric contravariant tensor fields on $M$. This leads to an interesting
link with the cohomology of $\Vect(M)$ with the operator coefficients,
namely, with values in the space $\Hom(\cS^k,\cS^{k-\ell})$.
Notice that this is the next natural case comparing with the
Gel'fand-Fuchs cohomology 
of $\Vect(M)$ with coefficients in the modules of tensor fields on $M$.

Examples of the $\Vect(M)$-cohomology classes
with values in $\Hom(\cS^k,\cS^{k-\ell})$ are given by the
1-cocycles $\g^\l_1$ and $\g^{1/2}_2$ (cf. Section \ref{Element}).
We will study this cohomology of $\Vect(M)$ in a subsequent article.

The restriction of the modules $\cD_\l$ to $\Sl(n+1,\bbR)$ leads to the
$\Sl(n+1,\bbR)$-cohomology (with the same coefficients). The complete answer
in this case was obtained in \cite{pierre}.

\bigskip
\noindent
{\it Acknowledgments}. It is a pleasure to acknowledge
numerous fruitful discussions
with C.~Duval and his constant interest to this work.
We are grateful to J.-L.~Brylinski, R.~Brylinski, M.~De~Wilde
and P.~Xu for
enlightening discussions.

%%%%%%%%%%%%%%%%%%%%%%%%%%%%%%%%%%%%%%%%%%%%%%%%%%%%%%%%%%%%%%%%%%%%%%%%%%%%%%

\end{document}